\documentclass[reqno,a4paper]{amsart}

\usepackage{amsmath,amsfonts,amssymb}
\usepackage{verbatim}
\usepackage{enumerate}
\usepackage{url}
\usepackage{tikz}
\usepackage{setspace}
\usepackage[left=2.5cm,right=2.5cm, top=3cm, bottom=2.5cm]{geometry}
\usepackage[utf8]{inputenc} 
\usepackage[T1]{fontenc}
\usepackage{hyperref}
\usetikzlibrary{arrows}

\def\Z{\mathbb Z}

\theoremstyle{plain}
\newtheorem{theorem}{Theorem}[section]
\newtheorem{lemma}[theorem]{Lemma}

\newtheorem{corollary}[theorem]{Corollary}
\newtheorem{proposition}[theorem]{Proposition}
\newtheorem{remark}[theorem]{Remark}

\newtheorem{conjecture}[theorem]{Conjecture}
\def\proof{\noindent {\it Proof: }}
\def\qed{\hfill\hbox{$\square$}}

\theoremstyle{definition}

\numberwithin{equation}{section}

\title{The $\{1,s\}$-weighted Davenport constant in $C_n^k$}

\author[F. E. Brochero Mart\'{\i}nez]{F. E. Brochero Mart\'{\i}nez}
\address{
Departamento de Matem\'{a}tica\\
Universidade Federal de Minas Gerais\\
UFMG\\
Belo Horizonte, MG\\
31270-901\\
Brazil\\
}
\email{fbrocher@mat.ufmg.br }

\author[S. Ribas]{S. Ribas}
\address{
Departamento de Matem\'{a}tica\\
Universidade Federal de Ouro Preto\\
UFOP\\
Ouro Preto, MG\\
35400-000\\
Brazil\\
}
\email{savio.ribas@ufop.edu.br}

\thanks{The first author was partially supported by FAPEMIG APQ-02973-17, Brazil.}

\date{\today}

\onehalfspace

\begin{document}

\maketitle

\begin{abstract}
Let $G$ be a finite abelian group and let $\varnothing \neq A \subset \Z$. The $A$-weighted Davenport constant of $G$ is the smallest positive integer ${\sf D}_A(G)$ such that every sequence $x_1 \boldsymbol{\cdot} {\dots} \boldsymbol{\cdot} x_{{\sf D}_A(G)}$ over $G$ has a non-empty subsequence $(x_{j_i})_i$ such that ${\varepsilon_1} x_{j_1} + {\varepsilon_2} x_{j_2} + {\dots} + {\varepsilon_t} x_{j_t} = 0$ for some $\varepsilon_1, \varepsilon_2, {\dots}, \varepsilon_t \in A$. In this paper, we obtain both upper and lower bounds for ${\sf D}_{\{1,s\}}(C_n^k)$, where $C_n$ denotes the cyclic group of order $n$, $s^2 \equiv 1 \pmod n$ and $s \not\equiv \pm1 \pmod n$. These bounds become sharp in some ``small'' cases.
\end{abstract}

\section{Introduction}

Given a finite abelian group $G$ written additively, the {\em zero-sum problems} study conditions to ensure that a given sequence over $G$ has a non-empty subsequence with prescribed properties (such as length, repetitions, weights) such that the sum of their terms equals $0$, the identity of $G$. This kind of problem dates back to the 60s, with the works of Erd\H os, Ginzburg \& Ziv \cite{EGZ} and Olson \cite{Ols1, Ols2}. It has been extensively studied for abelian groups, but recently several results over non-abelian groups have emerged. This paper deals with an weighted problem over $C_n^k$, where $C_n$ is the cyclic group of order $n$. First of all, we need some definitions and notations.

\subsection{Definitions and notations}

By a sequence $S$ over a finite group $G$ we mean a finite and unordered element of the free abelian monoid $\mathcal F(G)$ equipped with the sequence concatenation product denoted by $\boldsymbol{\cdot}$. A sequence $S \in \mathcal F(G)$ has the form $$S = \prod_{1 \le i \le k}^{\bullet} g_i = g_1 \boldsymbol{\cdot} {\dots} \boldsymbol{\cdot} g_k \in \mathcal F(G),$$ where $g_1, {\dots}, g_k \in G$ are the {\em terms} of $S$ and $k = |S| \ge 0$ is the {\em length} of $S$. Since the sequences are unordered, $S = \prod_{1 \le i \le k}^{\bullet} g_{\tau(i)} \in \mathcal F(G)$ for any permutation $\tau: \{1,2,{\dots},k\} \to \{1,2,{\dots},k\}$. For a given $S \in \mathcal F(G)$ and for $t \ge 0$, we abbreviate $S^{[t]} = \underbrace{S \boldsymbol{\cdot} {\dots} \boldsymbol{\cdot} S}_{t \text{ times}}.$
For $g \in G$, the {\em multiplicity} of the term $g$ in $S$ is denoted by $v_g(S) = \#\{i \in \{1,2,{\dots},k\} \, ; \; g_i = g\}$, therefore our sequence $S$ may also be written as $S = \prod_{g \in G}^{\bullet} g^{[v_g(S)]}.$
A sequence $T$ is a {\em subsequence} of $S$ if $v_g(T) \le v_g(S)$ for all $g \in G$; for that we use the notation $T \mid S$. In this case, we write $S \boldsymbol{\cdot} T^{[-1]} = \prod_{g \in G}^{\bullet} g^{[v_g(S) - v_g(T)]}$.
We also define:
\begin{align*}
\sigma(S) &= g_1 + {\dots} + g_k \in G, \quad \text{ the {\em sum} of $S$}; \\
\Sigma(S) &= \bigcup_{T \mid S \atop |T| \ge 1} \{\sigma(T)\} \subset G, \quad \text{ the {\em set of subsums} of $S$}; \\
S \cap K &= \displaystyle\prod_{{g \mid S \atop g \in K}}^{\bullet} g, \quad \text{ the subsequence of $S$ that lie in a subset $K$ of $G$}. 
\end{align*}
The sequence $S$ is called {\em zero-sum free} if $0 \not\in \Sigma(S)$, and $S$ is {\em zero-sum sequence} if $\sigma(S) = 0$.

\subsection{The Davenport constant}

One of the most important types of zero-sum problems is the {\em small Davenport constant} ${\sf d}(G)$ of a finite abelian group $G$. Namely, let ${\sf d}(G)$ be the maximal integer such that there exists a sequence over $G$ (repetition allowed) of length ${\sf d}(G)$ which is zero-sum free, i.e.,
$${\sf d}(G) = \sup\{|S|>0; \; S \in \mathcal F(G) \text{ is zero-sum free}\}.$$ 

Let $C_n$ be the cyclic group of order $n$. It is well-known that ${\sf d}(C_n) = n-1$. The inequality 
\begin{equation}\label{lowerboundabelian}
{\sf d}(C_{n_1} \oplus {\dots} \oplus C_{n_k}) \ge \sum_{i=1}^k (n_i - 1)
\end{equation}
holds true, considering the zero-sum free sequence formed by the concatenation of $n_i - 1$ copies of a generator of $C_{n_i}$ for each $1 \le i \le k$. Olson proved that the equality holds true for abelian $p$-groups \cite{Ols1} and for $k \le 2$ \cite{Ols2}.

Many variations and generalizations of zero-sum problems have been considered along the years. In this paper we consider the weighted problem: Let $A \subset \Z$ be a {\em set of weights}. In order to avoid trivial cases, we assume that $A \neq \varnothing$ and $A$ does not contain any multiple of $\exp(G)$, the exponent of $G$. Let
\begin{align*}
\sigma_A(S) &= \{{a_1} g_1 + {\dots} + {a_k} g_k \in G; \; a_i \in A\} \;\;\; \text{ the {\em set of $A$-weighted sums} of $S$}; \\
\Sigma_A(S) &= \bigcup_{T \mid S \atop |T| \ge 1} \sigma_A(T) \subset G \;\;\; \text{ the {\em set of $A$-weighted subsequence sums} of $S$}.
\end{align*}
Moreover, $S$ is {\em $A$-zero-sum free} if $0 \not\in \Sigma_A(S)$, and {\em $A$-zero-sum sequence} if $0 \in \sigma_A(S)$.

The {\em $A$-weighted Davenport constant} of an abelian group $G$ is defined by
$${\sf D}_A(G) = \inf\{k > 0; \; \text{ every $S \in \mathcal F(G)$ with $|S| \ge k$ is not $A$-zero-sum sequence}\}.$$ 

For instance, every abelian group $G$ yields ${\sf D}_{\{1\}}(G) = {\sf d}(G) + 1$, where ${\sf d}(G)$ is the (unweighted) small Davenport constant of $G$. Using Pigeonhole Principle, Adhikari {\it et al} \cite{Aea} proved that ${\sf D}_{\{\pm 1\}}(C_n) = \lfloor \log_2(n) \rfloor + 1$. More generally, Adhikari, Grynkiewicz \& Sun \cite{AGS} proved that if $n_1 \mid n_2 \mid {\dots} \mid n_k$ and $G = C_{n_1} \oplus {\dots} \oplus C_{n_k}$ then 
\begin{equation}\label{pesopm1dimd}
\sum_{i=1}^k \lfloor \log_2(n_i) \rfloor + 1 \le {\sf D}_{\{\pm 1\}}(G) \le \lfloor \log_2(|G|) \rfloor + 1,
\end{equation}
and Marchan, Ordaz \& Schmid \cite{MOS} removed the hypothesis $n_1 \mid n_2 \mid {\dots} \mid n_k$.

In addition, Adhikari, David \& Urroz \cite{ADU} showed that if $A = \{1, 2, 3, {\dots}, r\}$ then ${\sf D}_A(C_n) = \left\lceil \frac nr \right\rceil$ (it was previously proved for $n$ prime by Adhikari \& Rath \cite{AR}). Also in \cite{ADU}, it is proved that if $A$ is the set of all quadratic residues modulo a squarefree number $n$ and $\omega(n)$ is the number of distinct prime factors of $n$, then ${\sf D}_A(C_n) = 2 \omega(n) + 1$.

Although very little is known about the $A$-weighted Davenport constant for a general set $A$, Halter-Koch \cite{HK} gave an arithmetical interpretation of certain types of weighted Davenport constants in terms of algebraic integers and of binary quadratic forms.

In this paper, we fix an integer $k \ge 1$ and establish some upper and lower bounds for ${\sf D}_{\{1,s\}}(C_n^k)$, where 
\begin{equation}\label{hipotese}
s^2 \equiv 1 \!\!\!\pmod n, \quad \text{ but } \quad s \not\equiv \pm1 \!\!\!\pmod n.
\end{equation}
In fact, for $s \equiv 1 \pmod n$ we have ${\sf D}_{\{1,s\}}(C_n^k) = {\sf d}(C_n^k)+1$, and for $s \equiv -1 \pmod n$ we have ${\sf D}_{\{1,s\}}(C_n^k) = {\sf D}_{\{\pm1\}}(C_n^k)$. 

The main motivation for considering this set of weights was the inverse zero-sum problem related to small Davenport constant over the non-abelian group $C_n \rtimes_s C_2$. Indeed, in \cite{MR2}, the authors considered the inverse problem over $D_{2n}$, the dihedral group of order $2n$. The main argument used the $\{\pm1\}$-weighted Davenport constant over $C_n$ in order to obtain the structure of the product-one free sequences of maximum length. The idea for $C_n \rtimes_s C_2$ would be similar, using the $\{1,s\}$-weighted Davenport constant over $C_n$. However, the upper bounds provided in this paper were only able to solve the inverse problem for values of $n$ with sufficiently large factors, namely, factors $n_1$ and $n_2$ as in Lemma \ref{lemma} (in particular, the proof would not work for $n = 2^t$, $t \ge 3$). Nevertheless, in \cite{MR3} the authors completely solved the inverse problem over $C_n \rtimes_s C_2$, using the set of weights but without using directly the bounds on ${\sf D}_{\{1,s\}}(C_n^k)$ presented here.

This paper is organized as follows. In Section \ref{condicoes}, we present an useful factorization of $n$ in order to set a crucial projection, yielding an isomorphism in some cases. In Section \ref{lower}, we prove lower bounds for ${\sf D}_{\{1,s\}}(C_n^k)$. In Sections \ref{upper} and \ref{upper2}, we prove upper bounds for ${\sf D}_{\{1,s\}}(C_n^k)$ when the factorization of $n$ generates an isomorphism and a projection, respectively; in particular, we prove a relation among these two cases. In Section \ref{tight}, we discuss the tightness of these bounds and conclude the exact value in two families of ``small'' cases.

\section{The natural projection/isomorphism}\label{condicoes}

We note that conditions \eqref{hipotese} guarantee that $n$ can be neither an odd prime power nor twice an odd prime power, otherwise we would have $s \equiv \pm1 \pmod n$. The following lemma, which is also \cite[Lemma 2.2]{MR3} and we reproduce its proof here for convenience, ensures that these conditions suffice to factor $n$ nicely.

\begin{lemma}\label{lemma}
Let $n \ge 8$ and $s$ be positive integers satisfying the conditions \eqref{hipotese}.
\begin{itemize}
\item If both $n \neq p^t$ and $n \neq 2p^t$ for every prime $p$ and every integer $t \ge 1$, then there exist coprime integers $n_1, n_2 \ge 3$ such that $s \equiv -1 \pmod {n_1}$, $s \equiv 1 \pmod {n_2}$, and either (A) $n = n_1n_2$ or (B) $n = 2n_1n_2$.
\item If $n = 2^t$ for some $t \ge 3$, then (B) $n = 2n_1n_2$, where either $(n_1,n_2) = (1,2^{t-1})$ satisfies $s \equiv 1 \pmod {n_2}$ or $(n_1,n_2) = (2^{t-1},1)$ satisfies $s \equiv -1 \pmod {n_1}$.
\end{itemize}
\end{lemma}

\proof
Let $n = 2^tm$, where $m$ is odd and $t \ge 0$ is an integer. Since $m$ divides $s^2-1$ and $\gcd(s-1,s+1) \in \{1,2\}$ (depending on $s$ is even or odd), each prime power factor of $m$ divides either $s+1$ or $s-1$. Let $m_1 = \gcd(m,s+1)$ and $m_2 = \gcd(m,s-1)$, so that $m = m_1m_2$. In addition, $s^2 \equiv 1 \pmod {2^t}$ implies that either $s \equiv \pm1 \pmod {2^t}$ or $s \equiv 2^{t-1} \pm 1 \pmod {2^t}$ (whether $t \ge 3$). We consider some cases:
\begin{enumerate}[(i)]
\item {\bf CASE $t = 0$.} In this case, $n = m = m_1m_2$. We set $n_1 = m_1$ and $n_2 = m_2$, hence $n = n_1n_2$ is the desired factorization, as in (A). This is the only case where $s$ can be even; in the following, $s$ must be odd.
\item {\bf CASE $t = 1$.} It is possible to set either $n_1 = 2m_1$ and $n_2 = m_2$ or $n_1 = m_1$ and $n_2 = 2m_2$, hence $n = n_1n_2$ is the desired factorization, as in (A).
\item {\bf CASE $t \ge 2$ and $m \ge 3$.} If $s \equiv -1 \pmod {2^t}$, then we set $n_1 = 2^tm_1$ and $n_2 = m_2$. If $s \equiv 1 \pmod {2^t}$, then we set $n_1 = m_1$ and $n_2 = 2^tm_2$. Therefore, $n = n_1n_2$ is the desired factorization, as in (A). In the case that $t \ge 3$, it is possible that $s \equiv 2^{t-1} \pm 1 \pmod {2^t}$. If $s \equiv 2^{t-1}-1 \pmod {2^t}$, then we set $n_1 = 2^{t-1}m_1$ and $n_2 = m_2$. If $s \equiv 2^{t-1}+1 \pmod {2^t}$, then we set $n_1 = m_1$ and $n_2 = 2^{t-1}m_2$. Hence, $n = 2n_1n_2$ is the desired factorization as in (B).
\item {\bf CASE $t \ge 3$ and $m=1$.} In this case, $s \equiv 2^{t-1} \pm 1 \pmod {2^t}$, which implies that $s \equiv \pm1 \pmod {2^{t-1}}$. For the negative sign, it follows that $(n_1,n_2) = (2^{t-1},1)$, and for the positive sign, it follows that $(n_1,n_2) = (1,2^{t-1})$. Therefore, $n = 2 n_1 n_2$ is the factorization as in (B).
\end{enumerate}
\qed

By the previous lemma and Chinese Remainder Theorem, there exists a natural projection 
$$\Psi_0: C_n \to C_{n_1} \oplus C_{n_2}$$ 
satisfying 
$$\Psi_0(e) = (e_1,e_2) \quad \text{ and } \quad \Psi_0(s \cdot e) = (-e_1,e_2),$$ 
where $C_{n_1} = \langle e_1 \rangle$, $C_{n_2} = \langle e_2 \rangle$ and $C_n = \langle e \rangle$. In Case (A) of previous lemma, $\Psi_0$ is an isomorphism.

Let $\vec{0}_r$ be the identity of the group $C_r^k$, let $\vec{a} = (a_1,{\dots},a_k) \in C_r^k$, and write ${t \cdot \vec{a}} = (t a_1, {\dots}, t a_k)$. Denote by 
\begin{equation}\label{isomorphism}
\Psi: C_n^k \to C_{n_1}^k \oplus C_{n_2}^k 
\end{equation}
the similar projection satisfying 
$$\Psi(\vec{a}) = (\vec{a} , \vec{a}) \quad \text{ and } \quad \Psi(s \cdot \vec{a} ) = (-\vec{a} , \vec{a}),$$ 
where the first coordinate is the restriction to $C_{n_1}^k$ and the second coordinate is the restriction to $C_{n_2}^k$. In Case (A) of Lemma \ref{lemma}, $\Psi$ is an isomorphism. Moreover, for $g = (u,v) \in C_{n_1}^k \oplus C_{n_2}^k$, we write $s g = (- u,v)$.

When we look at the sequences over $C_{n_1}^k \oplus C_{n_2}^k$, since $s \equiv -1 \pmod {n_1}$ and $s \equiv 1 \pmod {n_2}$, in the second coordinate we have an unweighted problem, while in the first coordinate we have a $\{\pm1\}$-weighted problem. In the case that $n = n_1n_2$, finding zero-sum sequences over $C_n^k$ is equivalent to finding zero-sum sequences over $C_{n_1}^k \oplus C_{n_2}^k$. In the case that $n = 2n_1n_2$, zero-sums over $C_n^k$ imply zero-sums over $C_{n_1}^k \oplus C_{n_2}^k$, and $k+1$ disjoint zero-sums over $C_{n_1}^k \oplus C_{n_2}^k$ generate a zero-sum over $C_n^k$, since $C_n^k/(C_{n_1}^k \oplus C_{n_2}^k) \simeq C_2^k$ and ${\sf d}(C_2^k) = k$.

It is worth mentioning that the bounds presented here refer to the parities of $n_1$ and $n_2$ (see \cite[Conjecture~1.10]{GHST} for a conjecture regarding the Erd\H os-Ginzburg-Ziv constant over $C_n^3$ that is concerned with the parity of $n$). 
For the upper bound, our method consists basically of two steps: using the set of weights $\{\pm1\}$ in $C_{n_1}^k$ to produce more than ${\sf d}(C_{n_2}^k)$ disjoint sums in $\{\vec{0}_{n_1}\} \oplus C_{n_2}^k$, hence obtaining a zero-sum over $C_n^k$. However, the method can be improved provided $n_2$ is even in the following way: in the first step and with a little more effort, it is possible to obtain sums in $\{\vec{0}_{n_1}\} \oplus C_{\frac{n_2}{2}}^k$, where $C_{\frac{n_2}{2}} = \{2g; g \in C_{n_2}\}$. This effort is rewarded, since now it is only required ${\sf d}(C_{\frac{n_2}{2}}^k)+1$ sums for the second step. This method will be further applied for another zero-sum problems in future papers. 

\section{The lower bounds}\label{lower}

In this section, we provide good although not always tight lower bounds. The first one deals with the case $n = n_1n_2$.

\begin{theorem}\label{lowerthm}
Let $n$ and $s$ be as in \eqref{hipotese}, and write $n = n_1n_2$ as in Lemma \ref{lemma}. Then:
\begin{enumerate}[(1)]
\item $${\sf D}_{\{1,s\}}(C_n^k) \ge {\sf D}_{\{\pm1\}}(C_{n_1}^k) + {\sf d}(C_{n_2}^k) \ge k \left( \lfloor \log_2(n_1) \rfloor + n_2 - 1 \right) + 1;$$
\item If $n_2$ is odd and either $n_1 > n_2$ or $n_1$ is even, then $${\sf D}_{\{1,s\}}(C_n^k) \ge k(2n_2-1)+1.$$
\end{enumerate}
\end{theorem}

\proof
First of all, the projection defined in \eqref{isomorphism} ensures that we may consider $S$ over $C_{n_1}^k \oplus C_{n_2}^k$ instead of over $C_n^k$. For each case, we exhibit $\{1,s\}$-zero-sum free sequences $S$.

\begin{enumerate}[{(\it 1)}]
\item Let $S'_1$ be a $\{\pm1\}$-zero-sum free sequence over $C_{n_1}^k$ with $|S'_1| = {\sf D}_{\{\pm1\}}(C_{n_1}^k) - 1$ and $S'_2$ be a zero-sum free sequence over $C_{n_2}^k$ with $|S'_2| = {\sf d}(C_{n_2}^k)$. Set 
$$\quad \quad \quad S_1 = \prod_{g \in S'_1}^{\bullet} (g,\vec{0}_{n_2}) \in \mathcal F(C_{n_1}^k \oplus C_{n_2}^k), \quad S_2 = \prod_{h \in S'_2}^{\bullet} (\vec{0}_{n_1},h) \in \mathcal F(C_{n_1}^k \oplus C_{n_2}^k) \quad \text{and} \quad S = S_1 \boldsymbol{\cdot} S_2.$$
Then $S \in \mathcal F(C_{n_1}^k \oplus C_{n_2}^k)$ is $\{1,s\}$-zero-sum free and $|S| = |S_1| + |S_2| = {\sf D}_{\{\pm1\}}(C_{n_1}^k) + {\sf d}(C_{n_2}^k) - 1$, which proves (1).

\vspace{2mm}

\item Let $C_{n_1} = \left\langle e_1 \right\rangle$ and $C_{n_2} = \left\langle e_2 \right\rangle$. Set 
$$g_i = \Big( \underbrace{0,{\dots},0,e_1,0,{\dots},0}_{\text{$e_1$ in the $i^{\text{th}}$ position}} \; , \; \underbrace{0,{\dots},0,e_2,0,{\dots},0}_{\text{$e_2$ in the $(k+i)^{\text{th}}$ position}} \Big) \in C_{n_1}^k \oplus C_{n_2}^k \quad \text{for } 1 \le i \le k,$$
and 
$$S = \prod_{1 \le i \le k}^{\bullet} g_i^{[2n_2-1]} \in \mathcal F(C_{n_1}^k \oplus C_{n_2}^k).$$
Then $|S| = k(2n_2-1)$ and we claim that $S$ is $\{1,s\}$-zero-sum free. Suppose otherwise, then there exists a non-empty subsequence $T \mid S$ such that $0 \in \sigma_{\{1,s\}}(T)$. Assume that $g_i \mid  T$. By the $(k+i)^{th}$ position, we have $v_{g_i}(T) = n_2$. These $g_i$'s are the only possible terms of $T$ that modify the $i^{th}$ position. If $n_1$ is even then the $i^{th}$ position of the $\{1,s\}$-weighted sum is odd, hence it can not be $0$. If $n_1$ is odd and $n_1 > n_2$, then the $\{1,s\}$-weighted sum in the $i^{th}$ position is an odd element of $[-n_2,n_2]$, therefore it can not be $0$ modulo $n_1$.
\end{enumerate}
\qed

The second result deals with the case $n = 2n_1n_2$.

\begin{theorem}\label{lower2n1n2}
Let $n$ and $s$ be as in \eqref{hipotese}, and write $n = 2n_1n_2$ as in Lemma \ref{lemma}. Then ${\sf D}_{\{1,s\}}(C_n^k) \ge k \cdot {\sf D}_{\{1,s\}}(C_{n_1n_2}^k) + 1$. 
\end{theorem}

\proof
Let $T_1 \in \mathcal F(C_{n_1n_2})$ be a $\{1,s\}$-zero-sum free sequence of length $|T_1| = {\sf D}_{\{1,s\}}(C_{n_1n_2}) - 1$. The sequence $T_2 = T_1 \boldsymbol{\cdot} ((n_1n_2)e)$, where $C_{2n_1n_2} = \langle e \rangle$, is $\{1,s\}$-zero-sum free over $C_{2n_1n_2}$, and $|T_2| = {\sf D}_{\{1,s\}}(C_{2n_1n_2})$. Let $S \in \mathcal F(C_{2n_1n_2}^k)$ be the sequence formed by concatenation of the sequences $T_2$ in each coordinate, that is, 
$$S = \prod_{1 \le j \le k}^{\bullet} \prod_{g \mid T_2}^{\bullet} \underbrace{(0,\dots,0,g,0,\dots,0)}_{\text{$g$ in the $j^{th}$ coordinate}}.$$ 
Hence $|S| = k \cdot |T_2| = k \cdot {\sf D}_{\{1,s\}}(C_{n_1n_2})$ and $S$ is $\{1,s\}$-zero-sum free.
\qed

\section{The upper bounds for $n = n_1n_2$}\label{upper}

The projection $\Psi: C_n^k \to C_{n_1}^k \oplus C_{n_2}^k$ is actually an isomorphism in this case. Hence, finding zero-sum sequences over $C_n^k$ is equivalent to finding zero-sum sequences over $C_{n_1}^k \oplus C_{n_2}^k$. The goal in most cases in this section is to find many subsequences with a $\{\pm1\}$-zero-sum in the first coordinate, in order to obtain a zero-sum in the second coordinate.

Before claiming the main results of this section, we prove two simpler upper bounds.

\begin{proposition}\label{uppersimple}
Let $n$ and $s$ be as in \eqref{hipotese}, and write $n = n_1n_2$ as in Lemma \ref{lemma}. It holds:
\begin{enumerate}[(i)]
\item For every $n$, \quad $${\sf D}_{\{1,s\}}(C_n^k) \le \left({\sf d}(C_{n_2}^k) + 1\right) {\sf D}_{\{\pm1\}}(C_{n_1}^k) \le \left({\sf d}(C_{n_2}^k) + 1 \right) \Big( \lfloor k \log_2(n_1) \rfloor + 1 \Big);$$

\item For $n_2$ even, \quad $${\sf D}_{\{1,s\}}(C_n^k) \le \left({\sf d}(C_{\frac{n_2}2}^k) + 1\right) \Big( \lfloor k \log_2(n_1) \rfloor + k + 1 \Big).$$
\end{enumerate}
\end{proposition}

\proof
As in the proof of Theorem \ref{lowerthm}, we consider $S \in \mathcal F(C_{n_1} \oplus C_{n_2})$ instead of $S \in \mathcal F(C_n)$.

\vspace{1mm}

\begin{enumerate}[{\it (i)}]
\item 
Suppose that $|S| = \left({\sf d}(C_{n_2}^k) + 1 \right) {\sf D}_{\{\pm1\}}(C_{n_1}^k)$. Choose a minimal subsequence $A_1 \mid  S$ with $|A_1| \le {\sf D}_{\{1,s\}}(C_{n_1}^k)$ such that $\sigma_{\{1,s\}}(A_1) = (\vec{0}_{n_1}, b_1)$, where $b_1 \in C_{n_2}^k$. Inductively, we may construct $A_{j+1} \mid S \boldsymbol{\cdot} (A_1 \boldsymbol{\cdot} {\dots} \boldsymbol{\cdot} A_j)^{[-1]}$ for $j=1,{\dots}, {\sf d}(C_{n_2}^k)$ such that $|A_{j+1}| \le {\sf D}_{\{1,s\}}(C_{n_1})$ and $\sigma_{\{1,s\}}(A_{j+1}) = (\vec{0}_{n_1}, b_{j+1})$, where $b_{j+1} \in C_{n_2}^k$. On the other hand, it follows that we have ${\sf d}(C_{n_2}^k) + 1$ disjoint $\{1,s\}$-weighted sums in $\{\vec{0}_{n_1}\} \oplus C_{n_2}^k$, therefore we are done.

\vspace{1mm}

\item The proof of this inequality is similar to the previous one. The only difference is that each $\lfloor k \log_2(n_1) \rfloor + k + 1$ terms yield a $\{1,s\}$-weighted sum into $\{\vec{0}_{n_1}\} \oplus C_{\frac{n_2}{2}}^k$. In fact, let $T \mid S$ with $|T| = \lfloor k \log_2(n_1) \rfloor + k + 1$. We may construct $2^{|T|} > 2^k n_1^k$ sums using the terms of $T$, therefore the Pigeonhole Principle ensures that there exist distinct subsequences (but not necessarily pairwise disjoint) $T_1, T_2, {\dots}, T_{2^k+1} \mid T$ such that the first coordinates of $\sigma(T_1), \sigma(T_2), {\dots}, \sigma(T_{2^k+1})$ are all the same in $C_{n_1}^k$. Again by the Pigeonhole Principle, there exist two of them, say $T_1$ and $T_2$, such that $\sigma(T_1)$ and $\sigma(T_2)$ have the same parity in all their last $k$ coordinates. We may transform $T_1$ and $T_2$ into disjoint subsequences. Indeed, if $T_1$ and $T_2$ share the term $(a,b)$, then the first coordinate of $\sigma(T_1)$ equals the first coordinate of $\sigma(T_2)$ if and only if the first coordinate of $\sigma(T_1 \boldsymbol{\cdot} (a,b)^{[-1]})$ equals the first coordinate of $\sigma(T_2 \boldsymbol{\cdot} (a,b)^{[-1]})$, as well as all the last $k$ coordinates of $\sigma(T_1)$ have the same parity as those of $\sigma(T_2)$ if and only if all the last $k$ coordinates of $\sigma(T_1 \boldsymbol{\cdot} (a,b)^{[-1]})$ have the same parity as those of $\sigma(T_2 \boldsymbol{\cdot} (a,b)^{[-1]})$, where ``last $k$ coordinates'' denotes the restriction of the elements to $C_{n_2}^k$. Therefore, assuming that $T_1$ and $T_2$ are disjoint, we obtain $$\sum_{g \in T_1} sg + \sum_{h \in T_2} h \in \{\vec{0}_{n_1}\} \oplus C_{\frac{n_2}2}^k.$$
It is possible to construct ${\sf d}(C_{\frac{n_2}{2}}^k) + 1$ disjoint $\{1,s\}$-weighted sums as the previous one, therefore we are done.
\end{enumerate}
\qed

Our goal can be achieved faster once the subsequences with an $\{1,s\}$-zero-sum in the first coordinate have ``small'' length. For example, if we have two terms $(a_1,b_1) \boldsymbol{\cdot} (a_2,b_2) \mid S$ with $a_1 = \pm a_2$ then these two terms generate a $\{1,s\}$-zero-sum in the first coordinate. For three or four terms, we need the following:

\begin{lemma}\label{subseqle4}
Let $m \ge 2$ be integer, and let $S \in \mathcal F(C_m^k)$ with $|S| \ge \sqrt{2m^k}$. Suppose that every term of $S$ has multiplicity one, and the equation $g_1 + g_2 = 0$ has no solution with $g_1 \boldsymbol{\cdot} g_2 \mid S$. If $\vec{0}_m \nmid S$, then there exists $T \mid S$ with $3 \le |T| \le 4$ such that $\sum_{g \mid T} \varepsilon_g \cdot g = \vec{0}_m$ with $\varepsilon_g \in \{\pm1\}$.
\end{lemma}

\proof
Since $|S| \ge \sqrt{2m^k} > \sqrt{2m^k + 1/4} - 1/2$, we have that $\binom{|S|}{2} + \binom{|S|}{1} > m^k$. Therefore, the Pigeonhole Principle implies that $S$ contains distinct subsequences $T_1$ and $T_2$ of lengths one or two such that $\sigma(T_1) = \sigma(T_2)$. Notice that $T_1$ and $T_2$ have no common terms, otherwise either $T_1 = T_2$ or $\vec{0}_m \mid S$. Therefore, $\sigma(T_1) - \sigma(T_2) = \vec{0}_m$ and the latter has three or four terms.

\qed

\begin{theorem}\label{upperthm}
Let $n$ and $s$ be as in \eqref{hipotese}, and write $n = n_1n_2$ as in Lemma \ref{lemma}. It holds:
\begin{enumerate}[(I)]
\item For $n$ odd, \quad $${\sf D}_{\{1,s\}}(C_n^k) \le 2 {\sf d}(C_{n_2}^k) + \frac{n_1^k-1}{4} + \frac{\sqrt{2n_1^k}}{2};$$

\item For $n_1$ even, \quad $${\sf D}_{\{1,s\}}(C_n^k) \le 2 {\sf d}(C_{n_2}^k) + \frac{n_1^k + 2^k - 2}{4} + \frac{\sqrt{2n_1^k}}{2};$$

\item For $n_2$ even, \quad $${\sf D}_{\{1,s\}}(C_n^k) \le 2 {\sf d}(C_{\frac{n_2}2}^k) + \frac{n_1^k ( k \cdot 2^{k+1} + 1 ) + 2\sqrt{2n_1^k} + k (2^{k+1} + 4) + 7}{4(k+1)}.$$
\end{enumerate}
\end{theorem}

\proof
As in the proof of Theorem \ref{lowerthm}, we consider $S \in \mathcal F(C_{n_1} \oplus C_{n_2})$ instead of $S \in \mathcal F(C_n)$.

\begin{enumerate}[{\it (I)}]
\item 
In this case, $n_1$ and $n_2$ are both odd. Suppose that $|S| \ge 2 {\sf d}(C_{n_2}^k) + \frac{n_1^k-1}{4} + \frac{\sqrt{2n_1^k}}{2}$. \\
For $1 \le i \le u$, let $A_i = (\vec{0}_{n_1},b_i) \mid S$ be the terms of $S$ with $\vec{0}_{n_1}$ in the first coordinate. Furthermore, for $1 \le i \le v$, let $A_{u+i} = (a_{u+i}, b_{u+i}) \boldsymbol{\cdot} (a_{u+i}^{\pm1} , b'_{u+i}) \mid S$ denote the subsequences of $S$ formed by two terms such that their first coordinates are distinct than $\vec{0}_{n_1}$ and are equals or inverses in $C_{n_1}^k$. \\
From $S$, if we remove the terms of each $A_i$ with $1 \le i \le u+v$, it will remain at most $\frac{n_1^k-1}{2}$ terms (otherwise there would exist two terms equals or inverses in the first coordinate, that could be included in the previous $A_{u+i}$, $1 \le i \le v$). \\
Let $T = S \boldsymbol{\cdot} \left( \prod_{1 \le i \le u+v} A_i \right)^{[-1]}$, so that 
$$|T| = |S| - (u + 2v) \le \frac{n_1^k-1}{2}.$$ 
While $|T| \ge \sqrt{2n_1^k}$, Lemma \ref{subseqle4} ensures that there exists $T_1 \mid T$ with $|T_1| \in \{3,4\}$ such that $\vec{0}_{n_1} \in \sigma_{\{1,s\}}(T_1)$. Inductively, as long as possible, we construct $T_2, \dots, T_w$ satisfying the same conditions than $T_1$, hence $$|T| - 4w = |S| - (u + 2v + 4w) < \sqrt{2n_1^k}.$$
It follows that there exist $u+v+w$ disjoint subsequences of $S$ whose $\{1,s\}$-weighted sum belong to $\{\vec{0}_{n_1}\} \oplus C_{n_2}^k$. Since $$u + v + w = \frac{u+2v}{4} + \frac{u+2v+4w}{4} + \frac{u}{2} > \frac{|S| - \frac{n_1^k - 1}{2}}{4} + \frac{|S| - \sqrt{2n_1^k}}{4} \ge {\sf d}(C_{n_2}^k),$$
we are done.

\item 
This case is almost completely similar to the previous one. The only difference is the following: after defining $A_i$ for $1 \le i \le u+v$, the subsequence $T = S \boldsymbol{\cdot} \left( \prod_{1 \le i \le u+v} A_i \right)^{[-1]}$ satisfies 
$$|T| = |S| - (u + 2v) \le \frac{n_1^k + 2^k - 2}{2},$$ 
since there are $2^k-1$ non-zero elements of $C_{n_1}^k$ that are their own inverses. The remainder follows exactly the same steps than Case {\it(I)}.

\item 
In this case, $n_1$ is odd and $n_2$ is even. This proof is similar to the Case {\it (I)}. Suppose that $S \in \mathcal F(C_{n_1}^k \oplus C_{n_2}^k)$ with $|S| \ge 2 {\sf d}(C_{\frac{n_2}2}^k) + \frac{n_1^k ( k \cdot 2^{k+1} + 1 ) + 2\sqrt{2n_1^k} + k \cdot 2^{k+1} + 3}{4(k+1)}$. \\
For $1 \le i \le u$, let $A_i \mid S$ such that $|A_i| = 1$ and its only term lies in $\{\vec{0}_{n_1}\} \oplus C_{\frac{n_2}{2}}^k$. \\
For $1 \le i \le v$, let $A_{u+i} \mid S$ such that $|A_{u+i}| = 2$ and $\sigma_{\{1,s\}}(A_{u+i}) \in \{\vec{0}_{n_1}\} \oplus C_{\frac{n_2}{2}}^k$. Let $T_1 = S \boldsymbol{\cdot} \left( \prod_{1 \le i \le u+v} A_i \right)^{[-1]}$, so that 
$$|T_1| = |S| - (u+2v) \le 2^k \left(\frac{n_1^k - 1}{2}\right) + (2^k-1) = 2^{k-1}(n_1^k + 1) - 1.$$ 
For $1 \le i \le w$, let $A_{u+v+i} \mid T_1$ such that $|A_{u+v+i}| = 2$ and $\sigma_{\{1,s\}}(A_{u+v+i}) \in \{\vec{0}_{n_1}\} \oplus C_{n_2}^k$. Let $T_2 = T_1 \boldsymbol{\cdot} \left( \prod_{1 \le i \le u+v+w} A_i \right)^{[-1]}$, so that 
$$|T_2| = |S| - (u+2v+2w) \le \frac{n_1^k - 1}{2}.$$ 
For $1 \le i \le t$, let $A_{u+v+w+i} \mid S$ be as in Lemma \ref{subseqle4} in such way that $\sigma_{\{1,s\}}(A_{u+v+w+i}) \in \{\vec{0}_{n_1}\} \oplus C_{n_2}^k$. We have 
$$|S| - (u+2v+2w+4t) < \sqrt{2n_1^k}.$$ 
Since ${\sf d}(C_2^k) = k$, the disjoint subsequences $A_i$ for $1 \le i \le u+v+w+t$ produce at least $u+v+\left\lfloor \frac{w+t}{k+1} \right\rfloor$ terms into $\{\vec{0}_{n_1}\} \oplus C_{\frac{n_2}{2}}^k$. Since
\begin{align*}
\quad u+v+\frac{w+t}{k+1} &= \frac{u+2v+2w+4t}{4(k+1)} + \frac{u+2v+2w}{4(k+1)} + \frac{k(u+2v)}{2(k+1)} + \frac{u}{2} \\
&> \frac{|S| - \sqrt{2n_1^k}}{4(k+1)} + \frac{|S| - \frac{n_1^k-1}{2}}{4(k+1)} + \frac{k\left(|S| - 2^{k-1}(n_1^k + 1) + 1\right)}{2(k+1)} \ge {\sf d}(C_{\frac{n_2}{2}}^k) + 1,
\end{align*}
we are done.
\end{enumerate}
\qed

\begin{remark}
Notice that the previous theorem only used subsequences of length at most four in order to obtain a zero-sum in the first coordinate. It turns out that it is possible to go further with a similar of Lemma \ref{subseqle4}, considering subsequences of length at most six, eight, etc. However, the main terms involving ${\sf d}(C_{n_2}^k)$ and $n_1^k$ keep the same, the only change is that we get a larger denominator involving $\sqrt{2n_1^k}$ and we get other smaller terms of the form $\sqrt[3]{2n_1^k}$, $\sqrt[4]{2n_1^k}$, etc. For instance, if $n$ is odd and we consider subsequences up to six terms then we obtain the following inequality: $${\sf D}_{\{1,s\}}(C_n^k) \le 2 {\sf d}(C_{n_2}^k) + \frac{n_1^k - 1}{4} + \frac{\sqrt{2n_1^k}}{6} + \frac{\sqrt[3]{6n_1^k}}{3}.$$
Alternatively, we could use the remainder $\le \sqrt{2n_1^k}$ terms to form subsequences of length much larger, whose $\{1,s\}$-weighted sum is zero in the first coordinate and a better term (in some sense) in the second coordinate, but it also does not change the main terms. For these reasons, we decide to make it explicit until the first ``error term''.
\end{remark}

\section{The upper bounds for $n = 2n_1n_2$}\label{upper2}

Due to the projection $\Psi$, we have $C_n^k/(C_{n_1}^k \oplus C_{n_2}^k) \simeq C_2^k$. Therefore, one can simply multiply each bound of Theorem \ref{upperthm} by $k+1$, obtaining upper bounds for ${\sf D}_{\{1,s\}}(C_n^k)$. Indeed, it is possible to find $k+1$ disjoint subsequences whose $\{1,s\}$-weighted sums are $\vec{0}_{n_1} \oplus \vec{0}_{n_2}$, that is, belong to $C_2^k$. Since ${\sf d}(C_2^k) = k$, there exists a $\{1,s\}$-zero-sum subsubsequence. It proves the following result.

\begin{proposition}\label{prop2n1n2}
Let $n$ and $s$ be as in \eqref{hipotese}, and write $n = 2n_1n_2$ as in Lemma \ref{lemma}. Then 
$${\sf D}_{\{1,s\}}(C_n^k) \le (k+1) {\sf D}_{\{1,s\}}(C_{n_1n_2}).$$
\end{proposition}

In the case that $k=1$, $n_2$ is even and $s \equiv n_2 + 1 \pmod {2n_2}$, we are able to get a precise bound. For this, we need the following lemma:

\begin{lemma}[See {\cite[Theorem~11.1]{Gr1}}]\label{lemaciclico}
Let $m \ge 3$ be an integer and let $S \in \mathcal F(C_m)$ be a zero-sum free sequence with $|S| > m/2$. Then there exists $g \mid S$ such that $v_g(S) \ge \max\left\{ m - 2|S| + 1 , |S| - \left\lfloor \frac{m-1}{3} \right\rfloor \right\}$. In addition, for a residue class $r \pmod m$, denote by $\overline{r}$ the integer such that $0 \le \overline{r} \le m-1$ and $r \equiv \overline{r} \pmod m$. Then there exists an integer $t$ with $\gcd(t,m) = 1$ and $\sum_{g \mid S} \overline{gt} < m$. Moreover, for every $1 \le k \le \sum_{g \mid S} \overline{gt}$, there exists $T_k \mid S$ such that $\sum_{g \mid T_k} \overline{gt} = k$. 
\end{lemma}


We also need the value of ${\sf D}_{\{1,m+1\}}(C_{2m})$, where $m$ is even (so that $(m+1)^2 \equiv 1 \pmod {2m}$).

\begin{theorem}\label{s=m+1}
Let $m \ge 4$ be an even integer. Then 
$${\sf D}_{\{1,m+1\}}(C_{2m}) = m + 1.$$ 
Moreover, each sequence of length $m + 1$ has a proper subsequence whose $\{1,m+1\}$-weighted sum is $0 \in C_{2m}$.
\end{theorem}

\proof
The lower bound ${\sf D}_{\{1,m+1\}}(C_{2m}) \ge m + 1$ follows the same steps than Theorem \ref{lower2n1n2}, hence it is only required to prove the upper bound. Let $C_{2m} = \langle e \rangle$, so that $C_m = \langle 2e \rangle$. Let 
$$S = \prod_{1 \le i \le m+1}^{\bullet} (t_i e) \in \mathcal F(C_{2m}),$$ 
so that $|S| = m + 1$. We observe that ${\sf D}_{\{1,m+1\}}(C_m) = m = {\sf d}(C_m)+1$.

Since $|S| > m$, the projection $\Psi_1: C_{2m} \to C_m$ that maps $[(t \pmod {2m}) e]$ into $[(t \pmod m) e]$ implies that there exists $S_1 \mid S$ with $|S_1| \le m$ such that either $0 \in \sigma_{\{1,m+1\}}(S_1)$ or $me \in \sigma_{\{1,m+1\}}(S_1)$. If $0 \in \sigma_{\{1,m+1\}}(S_1)$, then we are done. Therefore we assume that $me \in \sigma_{\{1,m+1\}}(S_1)$. Notice that if $me \mid S$ then $|S \boldsymbol{\cdot} (me)^{[-1]}| = m$, and again we are done. Hence, we assume that $(me) \nmid S$. Furthermore, suppose that $S$ is zero-sum free over $C_{2m}$ (otherwise we are done). Since $|S| > |C_{2m}|/2$, after a possible change of generators (as in Lemma \ref{lemaciclico}), we assume that $t_1, \dots, t_{m+1}$ are positive integers in $[1,2m-1]$ such that $\sum_{i=1}^{m+1} t_i \le 2m-1$.

Notice that 
$$\begin{cases}
(m+1)a \equiv a \pmod {2m} &\text{ if $a$ is even,} \\
(m+1)a \equiv m + a \pmod {2m} &\text{ if $a$ is odd.}
\end{cases}$$
In this way, $S_1$ can be taken as an unweighted sum, that is, $\sigma(S_1) = me$. If $S_1$ has some odd term, say $ae \mid S_1$, then 
$$(m+1)a + \sum_{t_ie \mid (S_1 \boldsymbol{\cdot} (ae)^{[-1]})} t_i \equiv 0 \pmod {2m},$$ 
hence we are done. Otherwise, we assume that every term of any sequence $S_1$ with $\sigma(S_1) = me$ is even.

{\it CLAIM:} $(2e) \mid S_1$ and $e^{[2]} \mid S \boldsymbol{\cdot} S_1^{[-1]}$. \\
If this is the case, we consider $S_2 = (S_1 \boldsymbol{\cdot} e^{[2]}) \boldsymbol{\cdot} (2e)^{[-1]}$, therefore we are done. 

{\it PROOF OF THE CLAIM:} If $e^{[2]} \nmid S$ then 
$$t_1 + \dots + t_{m+1} \ge 1 + 2(m-1) = 1 + 2m > 2m,$$ 
a contradiction. If $2e \nmid S_1$ then $t_ie \mid S_1$ implies $t_i \ge 4$ is even. We have that $4 |S_1| \le \sum_{t_ie \mid S_1} t_i = m$, hence $|S_1| \le \frac{m}{4}$. Therefore 
$$|S \boldsymbol{\cdot} S_1^{[-1]}| = |S| - |S_1| \ge m + 1 - \frac{m}{4} = \frac{3m}{4} + 1,$$ 
which implies 
$$\sum_{t_ie \mid S \boldsymbol{\cdot} S_1^{[-1]}} t_i \ge \frac{3m}{4} + 1.$$ 
Let $v_1 = v_e(S) = v_e(S \boldsymbol{\cdot} S_1^{[-1]})$. Since $\sum_{t_ie \mid S \boldsymbol{\cdot} S_1^{[-1]}} t_i \le m - 1$ and $|S_1| \le \frac{m}{4}$, the average term of the subsequence $S \boldsymbol{\cdot} (S_1 \boldsymbol{\cdot} e^{[v_1]})^{[-1]}$ is at least $2$. Therefore
$$\frac{m-1-v_1}{\frac{3m}{4} + 1 - v_1} \ge \frac{\sum_{t_ie \mid S \boldsymbol{\cdot} S_1^{[-1]}} t_i - v_1}{|S \boldsymbol{\cdot} S_1^{[-1]}| - v_1} \ge 2 \Rightarrow v_1 \ge \frac{m}{2} + 3.$$
If $|S_1| = 1$ then $S_1 = me \mid S$, a contradiction. Therefore $|S_1| \ge 2$. The Pigeonhole Principle ensures that there exists $t_i e \mid S_1$ such that $1 \le t_i \le \frac{m}{2}$. In this case, we replace the term $t_ie$ by $t_i$ copies of $e$, which is possible since $v_1 > \frac{m}{2} \ge t_i$. This leads us to a contradiction, hence we are done.

%
\qed

\begin{corollary}\label{igual2n2par}
Let $n = 2n_1n_2$ and $s$ be as in Lemma \ref{lemma} such that $n_2$ is even and $s \equiv n_2+1 \pmod {2n_2}$. Then 
$${\sf D}_{\{1,s\}}(C_n) = n_1n_2 + 1.$$
\end{corollary}

\proof
In view of Theorem \ref{lower2n1n2}, it is only required to prove the upper bound. Let $S \in \mathcal F(C_n)$ such that $|S| = {\sf D}_{\{1,s\}}(C_{n_1n_2}) + 1 = n_1n_2 + 1$. Since $|S| \ge n_2+1$, Theorem \ref{s=m+1} ensures that there exists a subsequence $T_1 \mid S$ with $|T_1| \le n_2$ whose $\{1,s\}$-weighted sum belongs to the subgroup $C_{n_1}$. Now construct $T_2 \mid S \boldsymbol{\cdot} T_1^{[-1]}$ with the same property. Inductively, for $j \le n_1$, we construct $T_j \mid S \boldsymbol{\cdot} (T_1 \boldsymbol{\cdot} \dots \boldsymbol{\cdot} T_{j-1})^{[-1]}$ with $|T_j| \le n_2$ such that an $\{1,s\}$-weighted sum belongs to $C_{n_1}$. Since ${\sf d}(C_{n_1}) = n_1-1$, this corollary follows.

\qed


%
%
%
%
%
%
%

\section{The tightness of the bounds}\label{tight}

For the case $n = n_1n_2$, if $n_2$ is odd and either $n_1 > n_2$ or $n_1$ is even, then the second bound of Theorem \ref{lowerthm} is useful only when $n_2 \ge \lfloor \log_2(n_1) \rfloor$. Otherwise, the first bound is more useful. The worst bounds are those where $n_2$ is odd and $n_1 < n_2$, in which case we were unable to find a good lower bound. On the other hand, Proposition \ref{uppersimple} is more useful than Theorem \ref{upperthm} only when $n_1$ is very small compared to $n_2$; otherwise, Theorem \ref{upperthm} turns out to be of better use. 

For the case $n = 2n_1n_2$, we provided a good lower bound in terms of the previous case. In the particular case $k=1$, the problem splits into two subcases: $n_2$ even with $s \equiv n_2 + 1 \pmod {2n_2}$ (whose exact value has been established) and $n_1$ even with $s \equiv n_1 - 1 \pmod {2n_1}$. In both cases, the values of ${\sf D}_{\{1,s\}}(C_n)$ depend on ${\sf D}_{\{1,s\}}(C_{n_1n_2})$. On the other hand, the upper bound is asymptotically sharp for large values of $k$.

\vspace{1mm}

\noindent{\bf Notation.} Let $f, g$ be positive functions. As usual, we use the asymptotic Bachmann-Landau notations: 
\begin{itemize}
\item $f(n) \sim g(n)$ means that $\displaystyle\lim_{n \to \infty} \frac{f(n)}{g(n)} = 1$, 
\item $f(n) = O(g(n))$ means that there exists $c > 0$ such that $f(n) \le c \cdot g(n)$ for every positive integer $n$, and 
\item $f(n) = o(g(n))$ means that $\displaystyle\lim_{n\to \infty} \frac{f(n)}{g(n)} = 0$. 
\end{itemize}

\begin{conjecture}\label{conjec}
Let $n$ and $s$ be as in \eqref{hipotese}, and write $n = n_1n_2$ as in Lemma \ref{lemma}. We expect that 
\begin{align*}
{\sf D}_{\{1,s\}}(C_n^k) &= 
\begin{cases}
2 k n_2 + O(k \cdot \log_2(n_1)) \quad \text{ if $n_2$ is odd,} \\
k n_2 + O(k \cdot \log_2(n_1)) \quad\;\, \text{ if $n_2$ is even,}
\end{cases} \\
{\sf D}_{\{1,s\}}(C_n) &= 
\begin{cases}
2 n_2 + \lfloor \log_2(n_1) \rfloor \quad \quad\, \text{ if $n_2$ is odd,} \\
n_2 + \lfloor \log_2(n_1) \rfloor \quad \quad\;\,\, \text{ if $n_2$ is even.}
\end{cases}
\end{align*}
In the case that $n = 2n_1n_2$, we expect that 
$${\sf D}_{\{1,s\}}(C_n^k) = k \cdot {\sf D}_{\{1,s\}}(C_{n_1n_2}^k) + 1.$$
\end{conjecture}

Remember that ${\sf d}(C_n) = n-1$ and ${\sf d}(C_n^2) = 2n-2$. It is conjectured that if $k \ge 3$ then equality in \eqref{lowerboundabelian} holds for $C_n^k$, that is, ${\sf d}(C_n^k) = k(n-1)$. Girard \cite{Gir} showed that ${\sf d}(C_n^k) \sim kn$ as $n \to \infty$. As a consequence of the upper bounds and Girard's result, we have the following:

\begin{corollary}
Let $n$ and $s$ be as in \eqref{hipotese}, and write $n = n_1n_2$ as in Lemma \ref{lemma}.
\begin{enumerate}[(i)]
\item Fixed $n_1$, it follows that 
$${\sf D}_{\{1,s\}}(C_n^k) = 
\begin{cases}
2 k n_2 (1 + o(1)) \quad \text{ if $n_2$ is odd,} \\
k n_2 (1 + o(1)) \quad\;\, \text{ if $n_2$ is even.}
\end{cases}$$
\item Fixed $n_2$, it follows that $${\sf D}_{\{1,s\}}(C_n^k) = O(k \cdot \log_2(n_1)).$$
\end{enumerate}
\end{corollary}


Next proposition and corollary yield the exact values of ${\sf D}_{\{1,s\}}(C_n)$ when $n_1 \in \{3,5\}$ and $n_2$ is even for $n = n_1n_2$ and for $n = 2n_1n_2$, respectively, proving particular cases of Conjecture \ref{conjec}. It is worth mentioning that Theorems \ref{lowerthm}{\em (1)} and \ref{upperthm}{\em (III)} yield 
$$n_2 + \lfloor \log_2(n_1) \rfloor \le {\sf D}_{\{1,s\}}(C_n) \le n_2 + \lfloor \log_2(n_1) \rfloor + 1$$ 
provided $n = n_1n_2$, $n_2$ is even and $n_1 \in \{3,5\}$, and similar intervals for the cases where $n = 2n_1n_2$.

\begin{proposition}\label{n1=3}
Let $n$ and $s$ be as in \eqref{hipotese}, and write $n = n_1n_2$ as in Lemma \ref{lemma}, where $n_1 \in \{3,5\}$ and $n_2$ even. Then 
$${\sf D}_{\{1,s\}}(C_n) = n_2 + \lfloor \log_2(n_1) \rfloor.$$
\end{proposition}

\proof
The lower bound ${\sf D}_{\{1,s\}}(C_n) \ge n_2 + \lfloor \log_2(n_1) \rfloor$ is given by item (1) of Theorem \ref{lowerthm}. Therefore we just need to show that if $S \in \mathcal F(C_{n_1} \oplus C_{n_2})$ and $|S| = n_2 + \lfloor \log_2(n_1) \rfloor$ then $S$ is not $\{1,s\}$-zero-sum free. Let $C_{n_1} = \left\langle e_1 \right\rangle$, $C_{n_2} = \left\langle e_2 \right\rangle$, $C_{\frac{n_2}2} = \left\langle 2e_2 \right\rangle$, and $e_2+C_{\frac{n_2}{2}} = \{te_2 \in C_{n_2}; \text{ $t$ is odd}\}$.

Let $k \in \Z$ such that $|S \cap (\{0\} \oplus C_{n_2})| = n_2 - k$. Since ${\sf d}(C_{n_2}) = n_2-1$, if $k \le 0$, then we are done. Hence, assume $k \ge 1$. 

\begin{enumerate}[(a)]
\item {\bf CASE} $n_1 = 3$, $|S| = n_2+1$: It holds $|S \cap (\{e_1,2e_1\} \oplus C_{n_2})| = k + 1 \ge 2$. From these, it is possible to obtain a zero-sum in the first coordinate using any subsequence formed by two terms $(g_1,h_1) \boldsymbol{\cdot} (g_2,h_2) \mid S \cap (\{e_1,2e_1\} \oplus C_{n_2})$. In fact, if $g_1 = g_2$ then $(g_1 - g_2, h_1+h_2) = (0, h_1+h_2)$, and if $g_1 = e_1$ and $g_2 = 2e_1$ then $(g_1+g_2,h_1+h_2) = (0,h_1+h_2)$.

Let $S = T_1 \boldsymbol{\cdot} T_2 \boldsymbol{\cdot} T_3$, where 
\begin{align*}
T_1 &\mid S \cap (\{0\} \oplus C_{n_2}), \\
T_2 &\mid S \cap (\{e_1,2e_1\} \oplus C_{\frac{n_2}2}), \text{ and } \\
T_3 &\mid S \cap (\{e_1,2e_1\} \oplus (e_2 + C_{\frac{n_2}2})).
\end{align*}
Then it is possible to obtain at least
\begin{align*}
\left\lfloor \dfrac{|T_1|}{2} \right\rfloor + \left\lfloor \dfrac{|T_2|}{2} \right\rfloor + \left\lfloor \dfrac{|T_3|}{2} \right\rfloor &= \left\lfloor \dfrac{n_2 - k}{2} \right\rfloor + \left\lfloor \dfrac{|T_2|}{2} \right\rfloor + \left\lfloor \dfrac{k + 1 - |T_2|}{2} \right\rfloor \\
&\ge \dfrac{n_2 - k - 1}{2} + \dfrac{|T_2| - 1}{2} + \dfrac{k + 1 - |T_2| - 1}{2} \ge \dfrac{n_2}{2} - 1 
\end{align*}
elements in $\{0\} \oplus C_{\frac{n_2}2}$. Since ${\sf d}(C_{\frac{n_2}2}) = \frac{n_2}2 - 1$, if the latter inequality is strict then we are done. Otherwise, $n_2 - k$, $|T_2|$ and $k + 1 - |T_2|$ are all odd, hence $k$ is odd. Therefore, after removing the sums of those pairs, it remains one term from each subsequence $T_1$, $T_2$ and $T_3$. Let $(0, t_1 e_2) \mid T_1$, $(u e_1, 2t_2 e_2) \mid T_2$ and $(v e_1, (2t_3-1)e_2) \mid T_3$ these remainder terms. If $t_1$ is even then we have one more term in $\{0\} \oplus C_{\frac{n_2}2}$, and we are done. Hence, we suppose that $t_1$ is odd. In this case, 
$$\big( 0, (t_1 + 2t_2 + 2t_3 - 1)e_2 \big) \; \in \; \sigma_{\{1,s\}}\big( (0, t_1 e_2) \boldsymbol{\cdot} (u e_1, 2t_2 e_2) \boldsymbol{\cdot} (v e_1, (2t_3-1)e_2) \big),$$ 
and again we have one more element in $\{0\} \oplus C_{\frac{n_2}2}$. Therefore, we are done.

\item {\bf CASE} $n_1 = 5$, $|S| = n_2+2$: As in the previous case, we have $|S \cap (\{e_1,2e_1,3e_1,4e_1\} \oplus C_{n_2})| = k+2 \ge 3$. Define $S = T_1 \boldsymbol{\cdot} T_2 \boldsymbol{\cdot}  T_3 \boldsymbol{\cdot} T_4 \boldsymbol{\cdot} T_5$, where 
\begin{align*}
T_1 \mid S \cap (\{0\} \oplus C_{n_2}), \quad \quad
T_2 &\mid S \cap (\{e_1,4e_1\} \oplus C_{e_2+\frac{n_2}{2}}), \quad \quad \;\;
T_3 \mid S \cap (\{e_1,4e_1\} \oplus C_{\frac{n_2}{2}}), \\
T_4 &\mid S \cap (\{2e_1,3e_1\} \oplus C_{e_2 + \frac{n_2}{2}}), \quad \quad
T_5 \mid S \cap (\{2e_1,3e_1\} \oplus C_{\frac{n_2}{2}}). 
\end{align*}
Notice that the terms in each $T_i$ can be grouped in pairs in order to obtain products in $\{0\} \oplus C_{\frac{n_2}{2}}$. Since $n_2$ and $|S|$ are even, $|T_2| + |T_3| + |T_4| + |T_5| = k+2$ and 
$$\sum_{i=1}^5 \left\lfloor \dfrac{|T_i|}{2} \right\rfloor \ge \dfrac{n_2 - k - 1}{2} + \sum_{i=2}^5 \dfrac{|T_i| - 1}{2} + \frac{1}{2} = \frac{n_2}{2} - 1,$$
it is possible to obtain at least $\frac{n_2}{2} - 1$ disjoint subsums in $\{0\} \oplus C_{\frac{n_2}2}$. If the latter inequality is strict then we are done. Otherwise, it remains at least one term in four of the subsequences $T_i$ for $1 \le i \le 5$. We observe that one term from each $T_2$ and $T_3$ generates a subsum in $\{0\} \oplus C_{n_2}$, as well as one term from each $T_4$ and $T_5$ generates a subsum in $\{0\} \oplus C_{n_2}$. If $|T_1|$ is even then it remains the terms from $T_2$, $T_3$, $T_4$ and $T_5$, and we obtain two subsums in $\{0\} \oplus C_{n_2}$, which produce one more subsum in $\{0\} \oplus C_{\frac{n_2}{2}}$. Otherwise, either $T_2$ and $T_3$ generate one more subsum in $\{0\} \oplus C_{n_2}$ or $T_4$ and $T_5$ generate one more subsum in $\{0\} \oplus C_{n_2}$, and adding the remainder term from $T_1$ led us to one more subsum in $\{0\} \oplus C_{\frac{n_2}{2}}$. Hence we are done.
\end{enumerate}
\qed

\begin{remark}
The previous argument can be adapted to higher values of $n_1$, but the inequality obtained will be weaker. Thereby, it will be required to produce more subsums in $\{0\} \oplus C_{\frac{n_2}{2}}$ using at most one remainder term from each subsequence $T_1, \dots, T_{n_1}$. This creates several cases as $n_1$ grows.
\end{remark}

From previous proposition and Corollary \ref{igual2n2par}, it follows that:

\begin{corollary}
Let $n$ and $s$ be as in \eqref{hipotese}, and write $n = 2n_1n_2$ as in Lemma \ref{lemma}, where $n_1 \in \{3,5\}$ and $n_2$ even. Then 
$${\sf D}_{\{1,s\}}(C_n) = n_2 + \lfloor \log_2(n_1) \rfloor + 1.$$
\end{corollary}

\end{document}